\documentclass{amsart}
\usepackage{amssymb, amsthm}

\newtheorem{thm}{Theorem}[section]

\newtheorem{prop}[thm]{Proposition}

\theoremstyle{definition}
\newtheorem{defn}[thm]{Definition}
\newtheorem{exmp}[thm]{Example}

\theoremstyle{remark}
\newtheorem{rem}[thm]{Remark}

\def\CM{Cohen--Macaulay }
\def\der{\mathop{\mathrm{der}}}
\def\deg{\mathop{\mathrm{deg}}}
\def\div{\mathop{\mathrm{div}}}
\def\Proj{\mathop{\mathrm{Proj}}}
\def\Im{\mathop{\mathrm{Im}}}
\def\HH{\underline{\mathrm{Hom}}}

\def\H{\mathop{\mathrm{Hom}_R}}
\def\HS{\mathop{\mathrm{Hom}_S}}
\def\Cl{\mathop{\mathrm{Cl}}}
\def\zero{^{^{\circ}}}

\begin{document}

\title{Deformation of F-purity and F-regularity}

\author{Anurag K. Singh}

\address{Department Of Mathematics, The University Of Michigan, East Hall,
525 East University Avenue, Ann Arbor, MI 48109-1109}
\email{binku@math.lsa.umich.edu}


\date{November 21, 1997}

\begin{abstract}
For a Noetherian local domain $(R,m,K)$, it is an open question whether strong
F--regularity deforms. We provide an affirmative answer to this question when 
the canonical module satisfies certain additional assumptions. The techniques 
used here involve passing to a Gorenstein ring, using an anti--canonical cover. 
\end{abstract}

\maketitle

\section{Introduction}

Throughout our discussion, all rings are commutative, Noetherian, and have an 
identity element. Unless specified otherwise, we shall be working with rings 
containing a field of characteristic $p > 0$. The theory of the 
{\it tight closure}\/ of an ideal or a submodule of a module was developed by 
Melvin Hochster and Craig Huneke in \cite{HHjams} and draws attention to rings 
which have the property that all their ideals are tightly closed, called 
{\it weakly F--regular} rings. This property turn out to be of significant 
importance, for instance the Hochster--Roberts theorem of invariant theory that 
direct summands of polynomial rings are Cohen--Macaulay (\cite{HRinv}), can 
actually be proved for the much larger class of weakly F--regular rings. A 
closely related notion is that of a {\it strongly F--regular}\/ ring. The recent 
work of Lyubeznik and Smith (see \cite{LS}) shows that the properties of weak 
F--regularity and strong F--regularity agree for $\mathbb N$--graded F--finite 
rings.

\medskip

Whether the property of strong F--regularity deforms has been a persistent 
open question, specifically, if $R/xR$ is strongly F--regular for some 
nonzerodivisor $x\in m$, is $R$ strongly F--regular? A similar question may be 
posed when $(R,m,K)$ is a ring graded by the nonnegative integers, and $x\in m$ 
is a homogeneous nonzerodivisor. We examine these questions using the idea of 
passing to an anti--canonical cover $S=\oplus_{i \ge 0}I^{(i)}t^i$ (where $I$ 
represents the inverse of the canonical module in $\Cl(R)$) and get positive 
answers when the symbolic powers $I^{(i)}$ satisfy the Serre condition $S_3$ for 
all $i \ge 0$ and $S$ is Noetherian. It was earlier known that strong 
F--regularity deforms for Gorenstein rings, \cite[Theorem 4.2]{HHbasec}, and our 
results here may be considered an extension of those ideas. Another case when 
it is known that strong F--regularity deforms is in the case of 
$\mathbb Q$--Gorenstein rings, see \cite{AKM} and for the characteristic zero 
case, \cite{Smdeform}. It is interesting to note 
that all currently known results on the deformation of strong F--regularity 
require somewhat restrictive hypotheses which are weaker forms of the Gorenstein
hypothesis --- while $\mathbb Q$--Gorenstein is a natural weakening of the 
Gorenstein condition, the hypotheses required in our results are most obviously 
satisfied when $R$ is Gorenstein!

\medskip

Another notion ubiquitous in characteristic $p$ methods is that of F--purity
i.e., the property that for a ring $R$ of characteristic $p$, the Frobenius map 
remains injective upon tensoring with any $R$-module, see \cite{HRcohom}.
The corresponding question about the deformation of F--purity has a negative 
answer, the first examples being provided by R.~Fedder in \cite{Fed-pure}. 
However Fedder points out that his examples are less than satisfactory in two 
ways: firstly the rings are not integral domains, and secondly his 
arguments provide relevant examples only for finitely many choices of the 
characteristic $p$ of $K$. We point out various examples that have already been
studied in different contexts which overcome both these shortcomings, but now 
another shortcoming stands out --- although the rings $R$ are domains (which 
are not F--pure), the F--pure quotient rings $R/xR$ are not domains. 
The obvious question then arises: if $R/xR$ is an F--pure domain, is $R$ 
F--pure? An even more interesting question would be whether the property of
being an F--pure normal domain deforms. 

\medskip

In section \ref{dim2} we shall use cyclic covers to study graded rings of 
dimension two, and develop a criterion for the strong F--regularity of such 
rings when the characteristic of the field is zero or a  \lq\lq sufficiently 
large\rq\rq \ prime. By a graded ring we shall mean, unless stated otherwise, a 
ring graded by the
nonnegative integers, with $R_0=K$, a field. We shall denote by $m=R_+$, the
homogeneous maximal ideal of $R$, and by a system of parameters for $R$, we
shall mean a sequence of homogeneous elements of $R$ whose images form a system
of parameters for $R_m$. In specific examples involving homomorphic images of 
polynomial rings, lower case letters shall denote the images of the 
corresponding variables, the variables being denoted by the upper case letters.
We shall also implicitly assume our rings are homomorphic images of Gorenstein 
rings, and so possess canonical modules.  
 
\section{Frobenius closure and tight closure}

Let $R$ be a Noetherian ring of characteristic $p > 0$. The letter $e$ denotes 
a variable nonnegative integer, and $q$ its $e\,$th power, i.e., $q=p^e$.
We shall denote by $F$, the Frobenius endomorphism of $R$, and by $F^e$, its 
$e\,$th iteration, i.e., $F^e(r)=r^q$. For an ideal 
$I=(x_1, \dots, x_n) \subseteq R$, let $I^{[q]}=(x_1^q, \dots, x_n^q)$. Note 
that $F^e(I)R=  I^{[q]}$, where $q=p^e$, as always. Let $S$ denote the ring $R$
viewed as an $R$--algebra via $F^e$. Then $S\otimes_R\_$ is a covariant functor
from $R$--modules to $S$--modules. If we consider a map of free modules 
$R^n \to R^m$ given by the matrix $(r_{ij})$, applying $F^e$ we get a map 
$R^n \to R^m$ given by the matrix $(r_{ij}^q)$. For an $R$--module $M$, note 
that the $R$--module structure on $F^e(M)$ is $r'(r\otimes m) = r'r \otimes m$, 
and $r' \otimes rm = r'r^q \otimes m$. For $R$--modules $N \subseteq M$, we use 
$N_M^{[q]}$ to denote $\Im(F^e(N) \to F^e(M))$.

\medskip

For a reduced ring $R$ of characteristic $p > 0$, $R^{1/q}$ shall denote the
ring obtained by adjoining all $q\,$th roots of elements of $R$. The ring $R$
is said to be {\it F--finite}\/ if $R^{1/p}$ is module--finite over $R$. Note
that a finitely generated algebra $R$ over a field $K$ is F--finite if and only
if $K^{1/p}$ is finite over $K$. 

\medskip

We shall denote by $R\zero$ the complement of the union of the minimal 
primes of $R$. We say $I=(x_1,\dots,x_n) \subseteq R$ is a 
{\it parameter ideal}\/  if the images of $x_1,\dots,x_n$ form part of a system 
of parameters in $R_P$, for every prime ideal $P$ containing $I$.

\begin{defn} 
Let $R$ be a ring of characteristic $p$, and $I$ an ideal of $R$. For an
element $x$ of $R$, we say that $x \in I^F$, the {\it Frobenius closure }\/ of 
$I$, if there exists $q=p^e$ such that $x^q \in I^{[q]}$. 

\medskip

For $R$--modules $N \subseteq M$ and $u \in M$, we say that $u \in N_M^*$, the 
{\it tight closure }\/ of $N$ in $M$, if there exists $c \in R\zero$ such 
that $cu^q \in N_M^{[q]}$ for all $q=p^e \gg 0$. It is worth recording this 
when $M=R$, and $N=I$ is an ideal of $R$.
For an element $x$ of $R$, we say that $x \in I^*$, the {\it tight closure }\/ 
of $I$, if there exists $c \in R\zero$ such that $cx^q \in I^{[q]}$ for 
all $q=p^e \gg 0$. If $I^*=I$, we say that the ideal $I$ is {\it tightly closed}.
\end{defn}

It is easily verified that $I \subseteq I^F \subseteq I^*$. Furthermore, $I^*$ 
is always contained in the integral closure of $I$, and is frequently much 
smaller.

\begin{defn}
A ring $R$ is said to be {\it F--pure }\/ if the Frobenius homomorphism
$F: M \to F(M)$ is injective for all $R$--modules $M$.

\medskip

A ring $R$ is said to be {\it weakly F--regular}\/ if every ideal of $R$ is 
tightly closed, and is {\it F--regular}\/ if every localization is weakly 
F--regular. A weakly F--regular ring is F--pure.

\medskip

An F--finite ring $R$ is {\it strongly F--regular}\/ if for every element
$c \in R\zero$, there exists an integer $q=p^e$ such that the $R$--linear 
inclusion $R \to R^{1/q}$ sending $1$ to $c^{1/q}$ splits as a map of 
$R$--modules.

\medskip

$R$ is said to be {\it F--rational}\/ if every parameter ideal of $R$ is 
tightly closed. 
\end{defn}

The equivalence of the notions of strong F--regularity, F--regularity, and weak 
F--regularity, in general, is a persistent open problem. It is known that
these notions agree when the canonical module is a torsion element of the 
divisor class group, see \cite{Ma} and \cite{Williams}, and so when this 
hypothesis is satisfied, we switch freely between these notions. Furthermore,
Lyubeznik and Smith have shown that graded F--finite rings are weakly F--regular 
if and only if they are strongly F--regular, \cite[Corollary 4.3]{LS}.

\begin{thm}

\item $(1)$ \quad An F--finite regular domain of characteristic $p$ is strongly 
F--regular.

\item $(2)$ \quad  A strongly F--regular ring is F--regular, and an F--regular 
ring is weakly F--regular.

\item $(3)$ \quad Let $S$ be a strongly (weakly) F--regular ring. If $R$ is a 
subring of $S$ which is a direct summand of $S$ as an $R$--module, then $R$ is 
strongly (weakly) F--regular. 

\item $(4)$ \quad An F--rational ring $R$ is normal. If in addition, $R$ is 
assumed to be he homomorphic image of a Cohen--Macaulay  ring, then $R$ is 
Cohen--Macaulay.

\item $(5)$ \quad An F--rational Gorenstein ring is F--regular. If it is 
F--finite, then it is also strongly F--regular. 
\label{L:longlist}
\end{thm}

\begin{proof}
(1) follows from the fact that the Frobenius endomorphism over a regular local 
ring is flat. For all the assertions about strong F--regularity see 
\cite[Theorem 3.1]{HHstrong}. The results on F--rationality (4) and (5) are part 
of \cite[Theorem 4.2]{HHbasec}. 
\end{proof}

\section{F--purity does not deform}

We recall a useful result, \cite[Proposition 5.38]{HRcohom}.

\begin{prop}
Let $R=K[X_1, \dots, X_n]/I$ where $K$ is a field of characteristic $p$ and
$I$ is an ideal generated by square--free monomials in the variables 
$X_1, \dots, X_n$. Then $R$ is F--pure.
\end{prop}

In \cite{Fed-pure} Fedder constructs a family of examples for which F--purity 
does not deform. He points out that the examples have two shortcomings, firstly 
that they do not work for large primes, and secondly that the rings are not 
integral domains. We next review a special case of Fedder's examples and point 
out that it does serve as a counterexample to the deformation of F--purity for 
all choices of the prime characteristic $p$ of the field $K$. 

\begin{exmp}
Let $R=K[U,V,Y,Z]/(UV, UZ, Z(V-Y^2))$. Fedder's work shows that when the
characteristic of the field $K$ is $2$, then $R$ is not F--pure, although $R/yR$
is F--pure. We shall show here that this does not depend on the prime
characteristic $p$ of $K$, by showing that $y^3z^4 \notin I$ where 
$I = y^2(u^2-z^4)R$ but $(y^3z^4)^p \in I^{[p]}$. As for the assertion that
$R/yR$ is F--pure, note that $R/yR \cong K[U,V,Z]/(UV, UZ, ZV)$ and so F--purity
follows from the above proposition.

\medskip

In $R/I^{[p]}$ we have $(y^3z^4)^p = vy^{3p-2}z^{4p}= vy^{3p-2}u^{2p} =0$
and it only remains to show that $y^3z^4 \notin I$. Consider the grading on $R$
where $u,v,y$ and $z$ have weights $2$, $2$, $1$ and $1$ respectively. If 
$y^3z^4 \in I$, we have an equation in $K[U,V,Y,Z]$ of the form
$$
Y^3Z^4 = A(Y^2(U^2-Z^4))+B(UV)+C(UZ)+D(Z(V-Y^2)). 
$$
The polynomial $A$ is of degree $1$ and so must be a linear combination of
$Y$ and $Z$.  Examining this modulo $(U,V)$, we see that $A= -Y$, and so
$Y^3U^2=B(UV)+C(UZ)+D(Z(V-Y^2))$. Dividing through by $U$, we have
$Y^3U=B(V)+C(Z)+D'(Z(V-Y^2))$ where $D=UD'$, but then $Y^3U \in (V,Z)$, a
contradiction. 
\end{exmp}

We next use the idea of rational coefficient Weil divisors to construct a large
family of two dimensional normal domains for which F--purity does not deform, 
independent of the characteristic of the field $K$. These examples, as
presented, are for the graded case but local examples are obtained by the
obvious localization at the homogeneous maximal ideal. We recall some notation 
and results from \cite{De, Wadem, Wadim2}.
 
\begin{defn}
By a {\it rational coefficient Weil divisor}\/  on a normal projective variety
$X$, we mean a linear combination of codimension one irreducible subvarieties of
$X$, with coefficients in ${\mathbb Q}$.
\end{defn}
 
Let $K(X)$ denote the field of rational functions on $X$. Each $f \in K(X)$ 
gives us a Weil divisor $\div (f)$ by considering its zeros and poles with
appropriate multiplicity. For a divisor $D$, we shall say $D \ge 0$ if each 
coefficient of $D$ is nonnegative. Recall that for a Weil divisor $D$ on $X$, 
we have 
$$
H^0(X,{\mathcal{O}}_X(D)) = \{f \in K(X): \div (f) + D \ge 0 \}.
$$
If $D = \sum n_iV_i$ with $n_i \in {\mathbb Q}$ is a rational coefficient Weil 
divisor on $X$, we shall set $ [D]= \sum [n_i]V_i$, where $[n]$ denotes the 
greatest integer less than or equal to $n$, and define 
${\mathcal{O}}_X(D) ={\mathcal{O}}_X([D])$.  With this notation, Demazure's 
result is:
 
\begin{thm} 
Let $R=\oplus _{n \ge 0}R_n$ be a graded normal ring where $R_0 = K$ is a field.
Then $R$ can be described by a rational coefficient Weil divisor $D$ on 
$X = \Proj(R)$ satisfying the condition that $ND$ is an ample Cartier divisor
for some integer $N > 0$, in the form
$$
R = R(X,D) = \oplus _{n \ge 0} H^0(X,{\mathcal {O}}_X(nD))T^n \subseteq K(X)[T],
$$
where $T$ is a homogeneous element of degree one in the quotient field of $R$.
This identification preserves the grading, i.e., 
$R_n = H^0(X,{\mathcal {O}}_X(nD))T^n$.
\label{L:demazure}
\end{thm}

For a rational coefficient Weil divisor $D=\sum(p_i/q_i)V_i$ for $p_i$ and 
$q_i$ relatively prime integers with $q_i > 0$, we define
$D' = \sum((q_i-1)/q_i)V_i$ to be the {\it fractional part\/} of $D$. Note that
with this definition of $D'$ we have, for a positive integer $n$,  
$-[-nD] = [nD+D']$.
 
\medskip

Let $D$ be a rational coefficient Weil divisor on $X$ such that for some 
integer $N > 0$, $ND$ is an
ample Cartier divisors, and let $R= R(X,D)$ as above. One can get the following 
descriptions of the symbolic powers of the canonical module and their highest 
local cohomology modules, see \cite{Wadem, Wadim2}. If $X$ is a 
smooth projective variety of dimension $d$ with canonical
divisor $K_X$ and $\omega$ is the canonical module of $R$, we have
$$
[\omega^{(i)}]_n = H^0(X,{\mathcal O}_X(i(K_X+D')+nD))T^n,
$$
$$
[H^{d+1}_m(\omega^{(i)})]_n= H^d(X,{\mathcal O}_X(i(K_X+D')+nD))T^n. 
$$
It is easily seen that we can identify the action of the Frobenius on the
$n\,$th graded piece of the injective hull $E_R(K) = [H^{d+1}_m(\omega)]_n$ of 
$K$ (as a graded $R$--module) with the map
$$
H^d(X,{\mathcal O}_X(K_X+D'+nD))T^n\to H^d(X,{\mathcal O}_X(p(K_X+D'+nD)))T^{pn}.
$$
The action of the Frobenius on the socle of $E_R(K)$  
corresponds to the piece of this map in degree zero, i.e.,
$$
H^d(X,{\mathcal O}_X(K_X+D'))\to  H^d(X,{\mathcal O}_X(p(K_X+D'))).
$$
 
\begin{rem} In the above notation, if the ring $R$ is F--pure, then the
action of the Frobenius on $E_R(K)$ is injective and so
$H^d(X,{\mathcal O}_X(p(K_X+D')))$ must be nonzero.
\label{L:fpure}
\end{rem}

\begin{exmp}
We first fix positive integers $n$ and $k$ such that $1/n + 2/k < 1$. 
Take ${\mathbb{P}}^1 =\Proj K[X,Y]$, choose $k$ distinct nonzero elements 
$\alpha_1, \dots, \alpha_k$ in the field $K$, and construct the rational 
coefficient Weil divisor on ${\mathbb{P}}^1$
$$
D = (1/n)V(X-\alpha_1Y)+ \dots +(1/n)V(X-\alpha_kY).
$$
Next construct the graded ring $R=R({\mathbb{P}}^1,D)$ where 
$$
R = \oplus_{i \ge 0}R_i = \oplus_{i \ge 0}H^0({\mathbb{P}}^1,{\mathcal O}_{{\mathbb{P}}^1}(iD))Y^i,
$$ see [De]. Let
$\Delta_i = (X-\alpha_1Y) \dots (X-\alpha_{i-1}Y)(X-\alpha_{i+1}Y)\dots (X-\alpha_kY)$
for $1 \le i \le k$ and $\Delta = (X-\alpha_1Y) \dots (X-\alpha_kY)$.
Then $R$ is generated by the elements $Y,A_1, \dots,A_k$ where 
$A_i = \Delta_i XY^n /\Delta$ for $1 \le i \le k$. Now for $i \neq j$, 
$A_iA_j \in YR$, and it is easily verified that 
$$
R/YR \cong K[A_1, \dots,A_k]/(A_iA_j:i \neq j).
$$ 
Since $(A_iA_j:i \neq j)$ is an ideal generated by square--free monomials,
the ring $R/YR$ is F--pure. To see that $R$ is not F--pure, by Remark
\ref{L:fpure} it suffices to see that 
$H^1({\mathbb{P}}^1,{\mathcal O}_{{\mathbb{P}}^1}(p(K_{{\mathbb{P}}^1}+D'))) = 0$.
Note that by Serre duality, the vector space dual of this is 
$H^0({\mathbb{P}}^1,{\mathcal O}_{{\mathbb{P}}^1}((1-p)(K_{{\mathbb{P}}^1}+D')))$, 
and it certainly suffices to show that this is zero. Since we picked $n$ and 
$k$ satisfying $1/n + 2/k < 1$, we get that 
$\deg (1-p)(K_{{\mathbb{P}}^1}+D')=k(p-1)(1/n + 2/k - 1)$ is negative, and so 
$H^0({\mathbb{P}}^1,{\mathcal O}_{{\mathbb{P}}^1}((1-p)(K_{{\mathbb{P}}^1}+D')))=0$.

\medskip

It is easy to see that the ring $R$ is actually F--rational. The elements
$Y, \ A_1+ \dots +A_k$ form a homogeneous system of parameters for $R$, and 
since $R/YR$ is F--pure, the ideal $(Y, \ A_1+ \dots +A_k)R$ is its own 
Frobenius closure. The a--invariant of $R$ is $a(R) = -1$, and so $R$ is 
F--rational by \cite[Theorem 7.12]{HHjalg}\label{L:fpurenotdef}
\end{exmp}

\section{Cyclic covers and anti--canonical covers}

If $(R,m,K)$ is a local or graded normal ring with an ideal $I$ of pure 
height one satisfying $I^{(n)}\cong R$, by the {\it cyclic cover}\/ of $R$ with 
respect to $I$, we mean the ring 
$$
S = R\oplus I \oplus I^{(2)}\oplus \dots \oplus I^{(n-1)}.
$$
We use the identification $I^{(n)}= uR$ to give $S$ a natural ring
structure, i.e., $S$ is the ring 
$$
S=R[It,I^{(2)}t^2, \dots, I^{(i)}t^i, \dots]/(ut^n-1)
 =R[It, \dots, I^{(n-1)}t^{n-1}]/(ut^n-1).
$$
If the characteristic of $K$ is relatively prime to $n$, the inclusion of 
$R$ in $S$ is \'etale in codimension one. This holds because for a height one 
prime $P$ of $R$, $IR_P=vR_P$ and so $S_P=R_P[vt]/(rv^nt^n-1)$ where 
$r \in R_P$ is a unit. This is implicit in \cite{Wadim2}. Note that 
this also shows that $S$ is normal.

This construction is most interesting when $I=\omega$ is the canonical ideal of
$R$, and $\omega$ is an element of finite order in the divisor class group
$\Cl(R)$. In this case the cyclic cover of $R$ with respect to $\omega$ is 
usually called the {\it canonical cover}. (As the reader has already discovered,
we will sometimes speak of an ideal of pure height one as an element of
$\Cl(R)$, although we really mean the class it represents.)

If $R$ is of dimension $d$, we have the isomorphism
$$
\H(H_m^d(\omega^{(i)}),E_R) \cong \omega^{(1-i)}
$$
where $E_R$ is the injective hull of the residue field $K$, and consequently
$$
\HS(H_{m_S}^d(S),E_S) \cong \H(H_m^d(S),E_R) \cong S.
$$
$S$ is therefore Gorenstein whenever it is Cohen--Macaulay. The following is a
very useful result of Watanabe, \cite[Theorem 2.7]{Wadim2}. 

\begin{thm}
Let $(R,m) \to (S,n)$ be a finite local homomorphism of normal rings which is 
\'etale in codimension one. Then if $R$ is strongly F--regular (F--pure), so is 
$S$.
\label{L:etale}\end{thm}

An immediate consequence of this is the following theorem:

\begin{thm}
Let $(R,m,K)$ be a \CM normal ring of dimension two with a canonical ideal 
$\omega$ satisfying $\omega^{(n)}\cong R$ for $n$ relatively prime to the
characteristic of $K$. 
Then $R$ is (strongly) F--regular if and only the cyclic cover
$S= R\oplus \omega \oplus \dots \oplus \omega^{(n-1)}$ is F--rational.
\label{L:finitecover}
\end{thm}

\begin{proof}
Since $\omega$ is a torsion element of the divisor class group of $R$, the
results of \cite{Ma} or \cite{Williams} show that F--regularity is equivalent 
to strong F--regularity. If $R$ is F--regular, Theorem \ref{L:etale} then 
shows that $S$ is F--regular and so is certainly F--rational. For the converse, 
since $S$ is of dimension two our earlier discussion proves that $S$ is 
Gorenstein, and consequently is F--regular. The ring $R$, being a direct summand 
of $S$, is then F--regular. 
\end{proof}

We next consider the case when the canonical ideal $\omega$ of the \CM normal
ring $(R,m,K)$ is not necessarily of finite order in the divisor class group 
$\Cl(R)$. In this case, we construct the {\it anti--canonical cover}\/ $S$ by 
taking an ideal $I$ of pure height one, which is the inverse of $\omega$ in 
$\Cl(R)$ and forming the symbolic Rees ring 
$$
S=\oplus_{i \ge 0}I^{(i)}t^i.
$$
Note that symbolic Rees rings in general need not be Noetherian. However
there is a theorem of Watanabe which is very useful when the ring is indeed
Noetherian, \cite[Theorem 0.1]{Waantic}. 

\begin{thm}
Let $(R,m,K)$ be a normal ring and $I$ an ideal of pure height one, which is 
the inverse of $\omega$ in $\Cl(R)$. Let $S=\oplus_{i \ge 0}I^{(i)}t^i$ be the
anti--canonical cover as above. Then if $S$ is Noetherian, $R$ is strongly
F--regular (F--pure) if and only if $S$ is strongly F--regular (F--pure).
\label{L:antican}
\end{thm}

Note that if $S$ as above is Noetherian and Cohen--Macaulay, then it is 
Gorenstein. This can be inferred from a local cohomology calculation in 
\cite{Waantic}, and is also \cite[Theorem 4.8]{GHNV}.

\section{Deformation}

We now consider the question whether strong F--regularity deforms for local 
rings, i.e., if $R/xR$ is strongly F--regular for some nonzerodivisor $x\in m$, 
then is $R$ strongly F--regular? It is known that F--rationality deforms, 
\cite[Theorem 4.2 (h)]{HHbasec}, and our work makes use of this in an essential 
way.

\begin{thm}
Let $(R,m,K)$ be a normal local ring such that if $I$ is an ideal of pure height
one representing the inverse of $\omega$ in $\Cl(R)$, then the symbolic powers 
$I^{(i)}$ satisfy the Serre condition $S_3$ for all $i \ge 0$ and the 
anti--canonical cover $S=\oplus_{i \ge 0}I^{(i)}t^i$ is Noetherian. 

If $R/xR$ is strongly F--regular for some nonzerodivisor $x\in m$, then $R$ is 
also strongly  F--regular.
\label{L:fregdef}
\end{thm}

\begin{proof}
We may, if necessary, change $I$ to ensure that $xR$ is not one of its minimal 
primes. Since we have assumed that $I^{(i)}$ is $S_3$, the natural maps 
give us isomorphisms $(I/xI)^{(i)} \cong I^{(i)}/xI^{(i)}$. Hence 
$$
S/xS \cong R/xR \oplus (I/xI) \oplus (I/xI)^{(2)}\oplus \cdots
$$
serves as an anti--canonical cover for $R/xR$. Theorem \ref{L:antican} now shows 
that $S/xS$ is strongly F--regular, and so is also
Cohen--Macaulay. Hence $S$ is Cohen--Macaulay and therefore Gorenstein.
As F--rationality deforms and $S/xS$ is strongly F--regular, we get that $S$ is
F--rational. However $S$ is Gorenstein, and so is actually strongly F--regular. 
Finally, $R$ is a direct summand of $S$ and so is strongly F--regular.
\end{proof}

We next specialize to the case when $R$ has a canonical ideal which is a torsion
element of $\Cl(R)$. Note that the work of \cite{AKM} already shows that
F--regularity deforms in this setting, but we believe it is still interesting to
examine deformation from the point of view of cyclic covers. Recall that we do 
not need to 
distinguish between F--regularity and strong F--regularity in this case. Our
previous result needed the rather strong hypothesis that the symbolic powers
$I^{(i)}$ are $S_3$, but in this case it turns out that this hypothesis is
equivalent to the F--regularity of $R$, at least when the characteristic $p$
is \lq\lq sufficiently large.\rq\rq

\begin{thm}
Let $(R,m,K)$ be a normal local ring, where $K$ is a field of characteristic 
$p$, such that the canonical ideal $\omega$ of $R$ is of finite order in 
$\Cl(R)$, and $R/xR$ is F--regular for some nonzerodivisor $x\in m$.

If the symbolic powers $\omega^{(i)}$ satisfy the Serre condition $S_3$, then 
$R$ is F--regular. The converse is also true provided $p$ does not divide the 
order of $\omega$ in $\Cl(R)$.
\end{thm}

\begin{proof}
We choose an ideal $I$ of pure height one which gives an inverse for $\omega$ 
in $\Cl(R)$ such that $xR$ is not one of the minimal primes of $I$. Since 
$\omega$ is of finite order in $\Cl(R)$, so is its inverse $I$, and the 
anti--canonical cover $S=\oplus_{i \ge 0}I^{(i)}t^i$ is certainly Noetherian. 
The symbolic powers $I^{(i)}$ are isomorphic to some symbolic powers 
$\omega^{(j)}$, and consequently satisfy the condition $S_3$, by which we can 
conclude that $R$ is F--regular from the previous theorem. For the converse, 
since $R$ is F--regular, by Theorem \ref{L:finitecover} so is the cyclic cover 
$S = R\oplus I \oplus I^{(2)}\oplus \dots \oplus I^{(n-1)}$ where $n$ is the 
order of $I$ in $\Cl(R)$. Hence the summands $I^{(i)}$ are, in fact,
Cohen--Macaulay. 
\end{proof}

\section {Two dimensional graded rings} \label{dim2}

Let $(R,m,K)$ be a graded ring where $m$ is the homogeneous maximal ideal. We 
shall follow the notation of Goto and Watanabe in \cite{GW}. 
We recall that the highest local cohomology module $H_m^d(R)$ of $R$, where
$\dim R =d$ may be identified with $ \varinjlim R/(x_1^t,\dots x_d^t) $
where $x_1, \dots, x_d$ is a system of parameters for $R$ and the maps are
induced by multiplication by $x_1 \ldots x_d$. If $R$ is Cohen--Macaulay, 
these maps are injective. $H_m^d(R)$ has the natural structure of a graded 
$R$--module,  
where $\deg [r+(x_1^t,\dots x_d^t)] = \deg r - \sum_{i=1}^d x_i$, when $r$ and 
the $x_i$ are taken to be forms. With this grading on $H_m^d(R)$, we define
the $a$--$invariant$ of $R$ to be the highest integer $a$ such that 
$[H_m^d(R)]_a$ is nonzero, see \cite{GW}. 

\medskip

In general, for graded $R$--modules $M$ and $N$, we may define the graded 
$R$--module $\HH(M,N)$, where $[\HH(M,N)]_i$ is the abelian group consisting
of all {\it graded\/} R--linear homomorphisms from $M$ to $N(i)$ where the
convention for the grading shift is $[N(i)]_j=[N]_{i+j}$ for all $j \in \mathbb Z$.
This gives $\HH(M,N)$ a natural structure as a graded $R$--module. 
The injective hull of $K$ in the category of graded $R$--modules is 
$\underline E_R(K) = \HH(R,K)$. Consequently for all graded $R$--modules $M$,
we have $\HH(M,\underline E_R(K)) = \HH(M,K)$. With this 
notation, $\omega=\HH(H_m^d(R),K)$ is a {\it graded canonical module}\/ for 
$R$. Note that since $\omega=\HH(H_m^d(R),K)$, we have 
$\omega = \oplus_{i \ge -a}[\omega]_i$.

\medskip

If $R$ is Gorenstein, it is easy to see that $\omega \cong R(a)$, the
isomorphism now being a graded isomorphism.

\medskip

If $R$ has a graded canonical module $\omega$ satisfying $\omega^{(n)}= uR$,
set $k=(\deg u)/n$. We may construct a graded cyclic cover $S$ as
$$
S=R[\omega t, \dots, \omega^{(n-1)}t^{n-1}]/(ut^n-1)
$$
where we set $\deg t = -k$. The ring $S$ has a $(1/n)\mathbb Z$ grading, and it
can be verified that this is in fact a $(1/n)\mathbb Z_{\ge 0}$ grading. 
Consequently $S$ may be graded by the nonnegative integers, and results 
holding for such gradings do apply here. We next show that $a(S)= -k$. 

\medskip

Note that we have a graded isomorphism
$$
S \cong R \oplus \omega(k) \oplus \omega^{(2)}(2k) \oplus \dots \oplus
\omega^{(n-1)}(nk-k). 
$$
Using the fact that $\HH(H_m^d(\omega^{(i)}),E_R) \cong \omega^{(1-i)}$, we get 
the graded isomorphisms
\begin{align*}
\omega_S &\cong \HH(H_m^d(S),K) \\
&\cong \omega \oplus R(-k) \oplus \omega^{(-1)}(-2k) \oplus \dots \oplus
     \omega^{(2-n)})(k-nk) \\
&\cong [\omega(k) \oplus R \oplus \omega^{(-1)}(-k) \oplus \dots \oplus
     \omega^{(2-n)}(2k-nk)](-k) \\
&\cong [\omega(k) \oplus R \oplus \omega^{(n-1)}(nk-k) \oplus \dots \oplus
     \omega^{(2)}(2k)](-k) \\
&\cong S(-k).
\end{align*}
We can conclude from this that $a(S) = -k$. 

\medskip

We next recall a result of Fedder, \cite[Theorem 2.10]{Fed-2dim}. 

\medskip

Let $(R,m,K)$ be a graded ring over a perfect field $K$. 
Then $\mathfrak{S} \subseteq \der_KR$ is said to be {\it $D$--complete}\/ if 
for every element $a \in R$ with $D(a)=0$ for all $D \in \mathfrak{S}$, we have
$D(a)=0$ for all $D \in \der_KR$. 
If $\mathfrak S$ is a $D$--complete set of homogeneous derivations, we set 
$$
d_{\mathfrak S}(R) = \sup \{\deg(D):D \in \mathfrak{S} \}.
$$
Fedder's Theorem then is:

\begin{thm}
Let $(R,m,K)$ be a graded two dimensional normal ring over a perfect field $K$. 
Assume that the characteristic $p$ of $K$ satisfies $d_{\mathfrak{S} }(R) < p$
for some $D$--complete set $\mathfrak{S}$ of homogeneous derivations.
Then $R$ is F--rational if and only if it has a negative $a$--invariant.
\end{thm}

There is a corresponding result in characteristic zero, 
\cite[Theorem 3.6]{Fed-2dim}.  

\begin{thm}
Let $(R,m,K)$ be a graded two dimensional normal ring where $K$ is a field of 
characteristic zero. Then $R$ is F--rational if and only if it has a negative 
$a$--invariant.
\end{thm}

We are now in a position to combine this idea with the theory of cyclic covers
to get the following result:

\begin{thm}
Let $(R,m,K)$ be a graded two dimensional normal ring over a perfect field $K$
of characteristic $p$. Assume that $d_{\mathfrak{S} }(R) < p$ for some 
$D$--complete set $\mathfrak S$ of homogeneous derivations. If $n$ is the order 
of the graded canonical module $\omega$ in the graded divisor class group of 
$R$, also assume that $n$ and $p$ are relatively prime, and let 
$\omega^{(n)}=uR$. Then $R$ is F--regular if and only if $\deg u > 0$. 
\label{fedcyclic}
\end{thm}

\begin{proof}
Consider the graded cyclic cover 
$S=R[\omega t, \dots, \omega^{(n-1)}t^{n-1}]/(ut^n-1)$. We observed in Theorem
\ref{L:finitecover} that $R$ is F--regular if and only if $S$ is F--rational. 
By the above result of Fedder, this holds if and only is $a(S) < 0$. But 
$a(S)= -(\deg u)/n$ and so $R$ is F--regular if and only $\deg u > 0$. \qed
\end{proof}

The characteristic $0$ version follows similarly.

\begin{thm}
Let $(R,m,K)$ be a graded two dimensional normal ring over a field $K$ of
characteristic $0$. Let $n$ be the order of the graded canonical module 
$\omega$ in the graded divisor class group of $R$ and $\omega^{(n)}=uR$. 
Then $R$ is F--regular if and only if $\deg u > 0$.
\end{thm}

\begin{exmp}
Let $R=K[t,t^4x,t^4x^{-1},t^4(x+1)^{-1}]$  where $K$ is a field of 
characteristic $p$ which we assume to be a suitably large prime. This is a ring 
which is F--rational but not F--pure, see \cite{HHjalg, Wadim2}. By mapping a 
polynomial ring onto it, we may write $R$ as 
$$
R=K[T,U,V,W]/(T^8-UV, T^4(V-W)-VW, U(V-W)-T^4W).
$$
This is graded by setting the weights of $t,u,v$ and $w$ to be $1$, $4$, $4$ 
and $4$ 
respectively. We shall see that $R$ is not F--regular as an application of 
the above results. Note that with this grading we have $a(R) = -1$. A 
graded canonical module for $R$ is 
$\omega = 1/t^3(v,w)R$. It is easy to compute that 
$\omega^{(2)} = 1/t^6(v^2,vw,w^2)R$ and $\omega^{(3)} = 1/t^9(v^2-2vw+w^2)R$. 
Since $\deg [1/t^9(v^2-2vw+w^2)] = -1$, Theorem \ref{fedcyclic} shows that $R$ 
is not F--regular. 

\medskip

The cyclic cover $S$ of $R$ is isomorphic to $K[T,Y,Z]/(T^4+YZ^2-Y^2Z)$ where 
the inclusion $R \to S$ is obtained by extending the map 
$$
u \mapsto yz^2, \ \  v \mapsto y^3+yz^2-2y^2z, \ \ w \mapsto z^3+y^2z-2yz^2,  
\ \  t \mapsto t.
$$ 
The ring $S$ has a grading where the weights of $t,y$ and $z$ are $1$, $4/3$ 
and $4/3$ respectively. The $a$--invariant is $a(S)=1/3$ and corresponds to the 
fact that $S$ is not F--rational.
\end{exmp}

\begin{exmp}
Let $T=K[X,Y,Z]/(X^3-YZ(Y+Z))$ where $K$ is a field of characteristic $p$,
assumed to be a sufficiently large prime. Let $T$ have the obvious grading
where each variable has weight $1$. We let $R$ be the subring
$$
R=K[X,Y^3,Y^2Z,YZ^2,Z^3]/(X^3-YZ(Y+Z)).
$$ 
This again is an example from \cite{HHjalg}, see also \cite{Wadim2}. With the 
inherited grading, $a(R)=-1$ and so $R$ is
F--rational. A graded canonical module for $R$ is $\omega = (y,z)R$. It is
easily verified that $\omega^{(2)} = (y^2,yz,z^2)R$ and $\omega^{(3)} = R$.
Consequently $R$ is not F--regular. 

\medskip

$T$ is in fact the cyclic cover of $R$, under the natural inclusion. Of course, 
$a(T)=0$ and $T$ is not F--rational.
\end{exmp}

\section*{Acknowledgments}

The author wishes to thank Mel Hochster for many valuable discussions.

\end{document}